\documentclass{amsart}

\usepackage{amsthm}
\usepackage{amssymb}
\usepackage{graphics}

\newtheorem{theorem}{Theorem}[section]
\newtheorem*{proposition*}{Proposition}
\newtheorem{proposition}[theorem]{Proposition}
\newtheorem*{theorem*}{Theorem}
\newtheorem{lemma}[theorem]{Lemma}

\theoremstyle{definition}
\newtheorem{definition}[theorem]{Definition}

\theoremstyle{remark}
\newtheorem{remark}[theorem]{Remark}

\newcommand{\Crosscap}{\mathrm{Crosscap}}
\newcommand{\Betti}{\beta_1}

\newcommand{\ConstantEdgepaths}{\Gamma_{\mathrm{const}}}

\newcommand{\NumTangles}{N}
\newcommand{\angleb}[1]{\langle #1 \rangle} 
\newcommand{\circleb}[1]{\langle #1 \rangle^{\circ}}
\newcommand{\Edge}{e}
\newcommand{\Edgepath}{\gamma}
\newcommand{\EdgepathSystem}{\Gamma}
\newcommand{\BasicEdgepath}{\lambda}
\newcommand{\BasicEdgepathSystem}{\Lambda}

\newcommand{\Diagram}{\mathcal{D}}

\begin{document}

\title{Crosscap numbers of pretzel knots}

\author{Kazuhiro Ichihara}
\address{%
College of General Education, 
Osaka Sangyo University, 
3--1--1 Nakagaito, Daito, Osaka 574--8530, Japan }
\email{ichihara@las.osaka-sandai.ac.jp}

\author{%
    Shigeru Mizushima}
\address{%
        Department of Mathematical and Computing Sciences, 
        Tokyo Institute of Technology, 
        2--12--1 Ohokayama, Meguro, 
        Tokyo 152--8552, Japan}
\email{mizusima@is.titech.ac.jp}


\date{}

\begin{abstract}
The crosscap number of a knot in the 3-sphere is defined as 
the minimal first Betti number of non-orientable subsurfaces bounded 
by the knot. 
In this paper, we determine the crosscap numbers of a large class of the pretzel knots. 
The key ingredient to obtain the result is 
the algorithm of enumerating all essential surfaces 
for Montesinos knots developed by Hatcher and Oertel. 
\end{abstract}

\maketitle

\section{Introduction}

In the 3-sphere $S^3$, 
every knot bounds an orientable surface, called a \textit{Seifert surface}. 
The \textit{genus} $g(K)$ of a knot $K$ 
is then defined as the minimal genus of Seifert surfaces for $K$.
On the other hand, any knot in $S^3$ also bounds a non-orientable subsurface in $S^3$: 
For, at least one of the two checkerboard surfaces for a diagram of the knot must be non-orientable.
In view of this, similarly to the genus of a knot,
B.E. Clark defined the crosscap number of a knot as follows.
\begin{definition}[Clark,\cite{Cl}]
The \textit{crosscap number} $\Crosscap(K)$ of a knot $K$ is defined as 
the minimal first Betti number of non-orientable subsurfaces bounded by $K$ in $S^3$.
\end{definition}
\noindent For completeness, we define $\Crosscap(K)=0$ if and only if $K$ is the unknot.

In general, it is rather difficult to decide the crosscap number for a given knot. 
As far as the authors know, 
results on the determination of  crosscap numbers are only as follows: 
A knot has crosscap number one if and only if it is $(2,n)$-cabled, shown by Clark \cite{Cl}. 
The knot $7_4$ in the knot table has crosscap number three; $\Crosscap(7_4)=3$, 
proved by Murakami and Yasuhara \cite{MY}. 
This gives the first example showing that the inequality; 
$\Crosscap(K) \le 2 g(K) + 1$, presented by Clark, is sharp.  
In \cite{T01}, Teragaito determined the crosscap number of genus one knots. 
Next, in \cite{T04}, he gave a formula for the crosscap numbers of torus knots. 
Recently, for two-bridge knots, Hirasawa and Teragaito 
established a practical algorithm to calculate their crosscap numbers in \cite{HT}. 

In this paper, we determine the crosscap numbers of a large class of the pretzel knots. 
%
Let $K=P(p_1,p_2,\ldots,p_\NumTangles)$ 
be a non-trivial pretzel knot constructed from $\NumTangles$ rational tangles 
corresponding to $(1/p_1, 1/p_2,\ldots, 1/p_\NumTangles)$, 
where each $p_i$ denotes an integer and
\begin{align*}
\textrm{the integer $p_i$ is other than $0,\pm1$. } \tag{*}
\label{Eq:Not01Condition}
\end{align*}
It is easily seen that, on the tuple $(p_1,p_2,\ldots,p_\NumTangles)$, 
either of the following conditions must be satisfied for $K$ to become a knot (not a link); 
\begin{itemize}
\item[(a)] among $p_1, p_2, \ldots, p_\NumTangles$, exactly one of them is even and the others are odd, 
\item[(b)] $\NumTangles$ is odd and all of $p_1, p_2, \ldots, p_\NumTangles$ are odd.
\end{itemize}
Then, our result is stated as follows.
%
%
\begin{theorem}
\label{Thm:Main:Crosscap}
The crosscap number of a non-trivial pretzel knot $K=P(p_1,p_2,\ldots,p_\NumTangles)$ satisfying (*)
is determined as
\[
\Crosscap(K)= 
\left\{
\begin{array}{lcl}
\NumTangles-1 & \hspace{1cm} & \mbox{when } (\textrm{a}) \mbox{ is satisfied, }\\ 
\NumTangles & & 				\mbox{when } (\textrm{b}) \mbox{ is satisfied.}
\end{array}
\right.
\]
\end{theorem}

The key ingredient to prove the theorem is 
the algorithm of enumerating all essential surfaces 
for Montesinos knots, which are knots obtained 
by connecting a number of rational tangles in line, 
developed by Hatcher and Oertel in \cite{HO}. 
%
In the next section, we will collect basic notions, 
review their algorithm briefly, 
and give some supplements about the main theorem. 
The main theorem will be proved in Section \ref{Sec:Proof}.

The authors would like to thank the participants of 
the Workshop on Crosscap Number 
held at Shikotsuko, National Park Resort Villages, Hokkaido, Japan 
from 1st to 4th August, 2006, for helpful suggestions on the earlier drafts. 
Especially they are grateful to Gengyu Zhang, who organized the fruitful workshop.

\section{Preliminary}
\label{Sec:Preliminary}

In this section, we will prepare terminologies used in the rest of the paper, 
briefly review the algorithm of Hatcher and Oertel,
and give supplements to our main theorem.

\subsection{Definitions}

For a knot $K$, 
there is a connected, possibly non-orientable subsurface 
embedded in $S^3$ which has $K$ as the boundary. 
We call the intersection of the surface and the exterior of the knot 
a \textit{spanning surface} for $K$. 
It is a compact connected surface embedded in the exterior of the knot. 


An embedded surface $F$ in a knot exterior 
is called \textit{incompressible} if $F$ has no compression disks,
\textit{boundary incompressible} if $F$ has no boundary-compression disks,
and \textit{essential} if $F$ is incompressible and boundary-incompressible.
The boundary $\partial F$ of $F$
may intersect the meridian of $K$ more than once. 
The minimal number of such intersection points is called 
the \textit{number of sheets} of $F$, 
and is denoted by $\sharp s (F)$ or simply by $\sharp s$. 
In the following, the first Betti number and the Euler characteristic of the surface $F$ 
are denoted by $\Betti(F)$ and $\chi(F)$. 

\subsection{The algorithm of Hatcher and Oertel}

Here we give a simple review of the algorithm of Hatcher and Oertel,
mainly for pretzel knots.
Please refer \cite{HO} or \cite{IM} for details.

%
The algorithm assumes that the knot $K$ is a Montesinos knot,
has no integer tangles,
and has $3$ or more rational tangles.
For the Montesinos knot $K$, 
the algorithm 
enumerates all boundary slopes of essential surfaces in the exterior of $K$. 
In their algorithm, 
embedded surfaces are expressed by ``edgepath systems'' on 
a particular 1-dimensional cellular complex, called ``diagram'' $\Diagram$,
lying on the $u$-$v$ plane. 
An edgepath system $\EdgepathSystem$ is composed of ``edgepaths''
$\{ \Edgepath_1, \Edgepath_2, \ldots, \Edgepath_\NumTangles \}$, 
and each ``edgepath" $\Edgepath_i$ represents subsurfaces 
around the $i$-th tangle of $K$. 


The diagram $\Diagram$ is a graph on the $uv$-plane
lying in the region $-1\le u \le 1$.
Vertices of $\Diagram$
are classified into three types:
a vertex $\angleb{p/q}$,
a vertex $\circleb{p/q}$
and a vertex $\angleb{1/0}$,
where $p/q$ is an arbitrary irreducible fraction.
Their coordinates are $(u,v)=((|q|-1)/|q|,p/q)$,
$(1,p/q)$ and $(-1,0)$ respectively.
For $p/q$ and $r/s$,
if $|ps-qr|=1$,
vertices $\angleb{p/q}$ and $\angleb{r/s}$ are connected by an edge
denoted by $\angleb{p/q}$\,--\,$\angleb{r/s}$.
In particular, 
an edge of the form $\angleb{z}$\,--\,$\angleb{z+1}$ for some integer $z$
is called a {\em vertical edge}.
Besides,
there are edges of the form
$\angleb{p/q}$\,--\,$\circleb{p/q}$,
which are called {\em horizontal edges}.
%
%
An edge of the diagram $\Diagram$ is also called a {\em complete edge}.
In contrast,
a segment on a complete edge which is a proper subset of the complete edge as a set
is called a {\em partial edge}.

%
%

%
%
  \begin{figure}[htb]

   \begin{minipage}{160pt} 
   \begin{center}
    \begin{picture}(56,120) 
    \scalebox{0.8}{
     \put(0,0){\scalebox{0.7}{\includegraphics{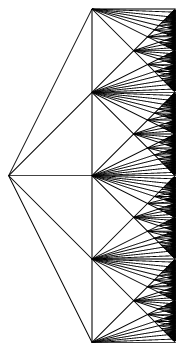}}}
     \put(49,135){\rotatebox{90}{\scalebox{1.0}{$\cdots$}}}
     \put(49,-10){\rotatebox{90}{\scalebox{1.0}{$\cdots$}}}
     \put(2,70.8){\scalebox{1.5}{\vector(1,0){55}}}
     \put(42,5){\scalebox{1.5}{\vector(0,1){90}}}
     \put(29,18){\tiny $-2$}
     \put(29,43){\tiny $-1$}
     \put(36,98){\tiny $1$}
     \put(36,123){\tiny $2$}
     \put(40,148){$v$}
     \put(3,82){\vector(1,-1){10}}
     \put(-10,85){$\angleb{1/0}$} 
     \put(2,64){\tiny $-1$}
     \put(72,64){\tiny $1$}
     \put(92,68){$u$}
     \put(33,63){\tiny $O$}
    }
    \end{picture}
   \end{center}
   \caption{The diagram $\Diagram$}
   \label{Fig:Diagram}
   \end{minipage}
    \begin{minipage}{160pt} 
    \begin{center}    
     \begin{picture}(120,120) 
     \scalebox{0.8}{
      \put(0,-10){\scalebox{0.75}{\includegraphics{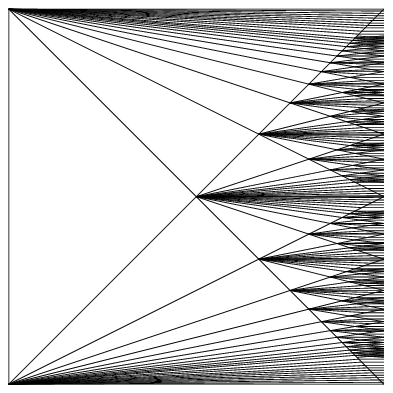}}}
      \put(-5,-2){\vector(1,0){10}}
      \put(-21,-5){$\angleb{0}$}
      \put(156,-2){\vector(-1,0){10}}
      \put(158,-5){$\circleb{0}$}
      \put(-5,133){\vector(1,0){10}}
      \put(-21,130){$\angleb{1}$}
      \put(156,133){\vector(-1,0){10}}
      \put(158,130){$\circleb{1}$}
      \put(63,66){\vector(1,0){10}}
      \put(46,63){$\angleb{\frac{1}{2}}$}
      \put(156,66){\vector(-1,0){10}}
      \put(158,63){$\circleb{\frac{1}{2}}$}
      \put(83,43){\vector(1,0){10}}
      \put(66,40){$\angleb{\frac{1}{3}}$}
      \put(156,43){\vector(-1,0){10}}
      \put(158,40){$\circleb{\frac{1}{3}}$}
      \put(94,32){\vector(1,0){10}}
      \put(77,29){$\angleb{\frac{1}{4}}$}
      \put(176,32){\vector(-1,0){30}}
      \put(178,29){$\circleb{\frac{1}{4}}$}
     }
     \end{picture}
    \end{center}
    \caption{
     A part of the diagram $\Diagram$ in 
     $[0,1]\times[0,1]$
    }
    \label{Fig:Diagram2}
    \end{minipage}

  \end{figure}

In the algorithm,
for a fixed Montesinos knot,
we first enumerate all ``basic edgepath systems''.
A {\em basic edgepath system} is a collection $\{\BasicEdgepath_1, \BasicEdgepath_2, \ldots, \BasicEdgepath_\NumTangles\}$ of $\NumTangles$ basic edgepaths.
As for a pretzel knot satisfying (\ref{Eq:Not01Condition}),
each {\em basic edgepath} $\BasicEdgepath_i$ is either of
$\BasicEdgepath_{i,a}=$
$\angleb{0}$\,--\,$\angleb{1/p_i}$
or
$\BasicEdgepath_{i,b}=$
$\angleb{s_i}$\,--\,$\angleb{s_i/2}$\,--\,$\ldots$\,--\,$\angleb{s_i/(|p_i|-1)}$\,--\,$\angleb{s_i/|p_i|}$,
where $s_i$ denotes the sign $+1$ or $-1$ of $p_i$.
Basic edgepaths $\BasicEdgepath_{i,a}$ and $\BasicEdgepath_{i,b}$ are 
extended to
$\widetilde{\BasicEdgepath_{i,a}}=$
$\angleb{0}$\,--\,$\angleb{1/p_i}$\,--\,$\circleb{1/p_i}$,
and
$\widetilde{\BasicEdgepath_{i,b}}=$
$\angleb{s_i}$\,--\,$\angleb{s_i/2}$\,--\,$\ldots$\,--\,$\angleb{s_i/(|p_i|-1)}$\,--\,$\angleb{s_i/|p_i|}$\,--\,$\circleb{s_i/|p_i|}$
which are called {\em extended basic edgepaths}.
In the same way,
we can extend a basic edgepath system 
$\BasicEdgepathSystem=\{ \BasicEdgepath_1, \BasicEdgepath_2, \ldots, \BasicEdgepath_\NumTangles \}$
to an extended basic edgepath system 
$\widetilde{\BasicEdgepathSystem}=\{ \widetilde{\BasicEdgepath_1}, \widetilde{\BasicEdgepath_2}, \ldots, \widetilde{\BasicEdgepath_\NumTangles} \}$.
An extended basic edgepath $\widetilde{\BasicEdgepath}$ can be regarded as a function $[0,1]\rightarrow\mathbb{R}$
which returns $v_0$ for $u_0$
when $\widetilde{\BasicEdgepath}$ meets with the vertical line $u=u_0$ 
at $(u_0,v_0)$,
and moreover,
an extended basic edgepath system
can be regarded as a function $[0,1]\rightarrow \mathbb{R}$
defined by $\widetilde{\BasicEdgepathSystem}(u)=\sum_{i=1}^\NumTangles \widetilde{\BasicEdgepath_i}(u)$.

Then,
next in the algorithm,
we enumerate {\em candidate edgepath systems} of type I, II or III.
(I)
For each basic edgepath system $\BasicEdgepathSystem$,
we solve the equation $\widetilde{\BasicEdgepathSystem}(u)=0$.
Then, for each solution $u_0$,
we make an edgepath system $\EdgepathSystem=\{\Edgepath_i\}$ as follows.
If $u_0 \le (|p_i|-1)/|p_i|$,
then we set $\Edgepath_i=\BasicEdgepath_i \cap \{(u,v)|u\ge u_0\}$,
where a partial edge may be included in $\Edgepath_i$.
If otherwise,
then we set $\Edgepath_i=\{P_i\}$
where $P_i$ is a point with $uv$-coordinate $(u,v)=(u_0,1/p_i)$
lying on a horizontal edge $\angleb{1/p_i}$\,--\,$\circleb{1/p_i}$.
This edgepath is called a {\em constant edgepath}.
An edgepath system thus defined is called 
a {\em type I edgepath system}.
(II)
For a basic edgepath $\BasicEdgepath$,
by adding vertical edges to $\BasicEdgepath$,
we can make an edgepath,
for instance,
$\angleb{-2}$\,--\,$\angleb{-1}$\,--\,$\angleb{0}$\,--\,$\angleb{1/p_i}$
or
$\angleb{s_i+1}$\,--\,$\angleb{s_i}$\,--\,$\angleb{s_i/2}$\,--\,$\ldots$\,--\,$\angleb{s_i/(|p_i|-1)}$\,--\,$\angleb{s_i/|p_i|}$.
Such an edgepath 
and a basic edgepath itself are called {\em type II edgepaths}.
For each basic edgepath system $\BasicEdgepathSystem$,
we make an edgepath system $\EdgepathSystem$
that 
each edgepath $\Edgepath_i$ of $\EdgepathSystem$
is a type II edgepath for $1/p_i$
and that $v$-coordinates of endpoints of the edgepaths in $\EdgepathSystem$ sum up to $0$.
Such an edgepath system is called a {\em type II edgepath system}.
(III)
For a basic edgepath $\BasicEdgepath_i$,
we make an edgepath
$\Edgepath_i=\angleb{1/0}$\,--\,$\BasicEdgepath_i$,
which is called a {\em type III edgepath}.
For each basic edgepath system $\BasicEdgepathSystem=\{\BasicEdgepath_i\}$,
we make an edgepath system $\EdgepathSystem$
that
each edgepath $\Edgepath_i$ of $\EdgepathSystem$
is a type III edgepath for $\BasicEdgepath_i$.
Such an edgepath system is called a {\em type III edgepath system}.

A candidate edgepath system corresponds to some properly embedded surfaces
in the exterior of $K$ as follows.
Now, we divide $S^3$ into $\NumTangles$ $3$-balls
so that each ball $B_i$ includes the $i$-th tangle.
Each edgepath $\Edgepath_i$ represents a subsurface $F_i$ in the ball $B_i$
whose boundary appears at the tangle in $B_i$ and on the boundary $\partial B_i$.
An edge in an edgepath expresses saddles.
A non-constant edgepath corresponds to subsurfaces 
obtained by combining saddles corresponding to the edges in $\Edgepath_i$.
Note here that two possible choices of saddles exist for one edge,
and that multiple subsurfaces correspond to an edgepath in general.
Besides,
a constant edgepath represents a {\em cap subsurface}.
Any candidate edgepath system $\EdgepathSystem$ satisfies conditions that
for the endpoints of edgepaths in the edgepath system,
$u$-coordinates coincide
and $v$-coordinates sum up to $0$.
These conditions assure that subsurfaces $\{F_i\}$ corresponding to edgepaths $\{\Edgepath_i\}$ are glued consistently
into a properly embedded surface $F$,
which is called a {\em candidate surface}.
These conditions are called the condition for {\em gluing consistency}.

Any essential surface can be isotopically deformed into some standard form,
and hence, is isotopic to one of the candidate surfaces
(In a precise sense,
some essential surfaces are represented by some exceptional edgepath system
not described above.
However, we can ignore such essential surfaces in this paper.)
Conversely, the set of candidate surfaces includes inessential surfaces in general.
Some conditions for an edgepath system to correspond to an essential surface
are given in \cite{HO}.
With these conditions,
by eliminating candidate edgepath systems which give only inessential surfaces,
we have the list of edgepath systems corresponding to essential surfaces.
By calculating boundary slopes for the essential surfaces,
we eventually obtain the list of boundary slopes.

\subsection{Supplements to the main theorem}

\begin{figure}[hbt]
 \begin{center}
  \begin{picture}(275,130)
   \put(0,18){\scalebox{0.3}{\includegraphics{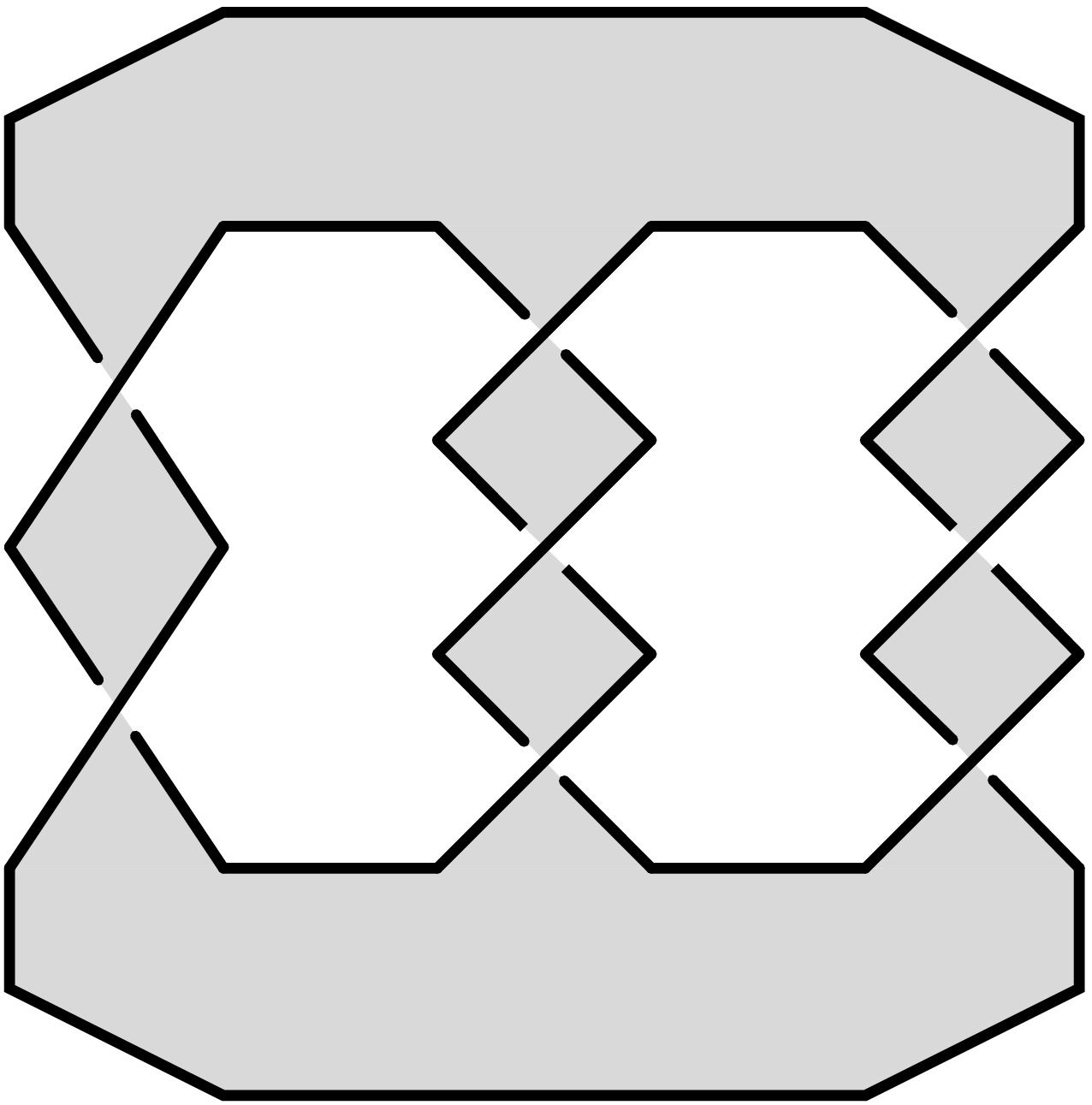}}}
   \put(160,0){\scalebox{0.3}{\includegraphics{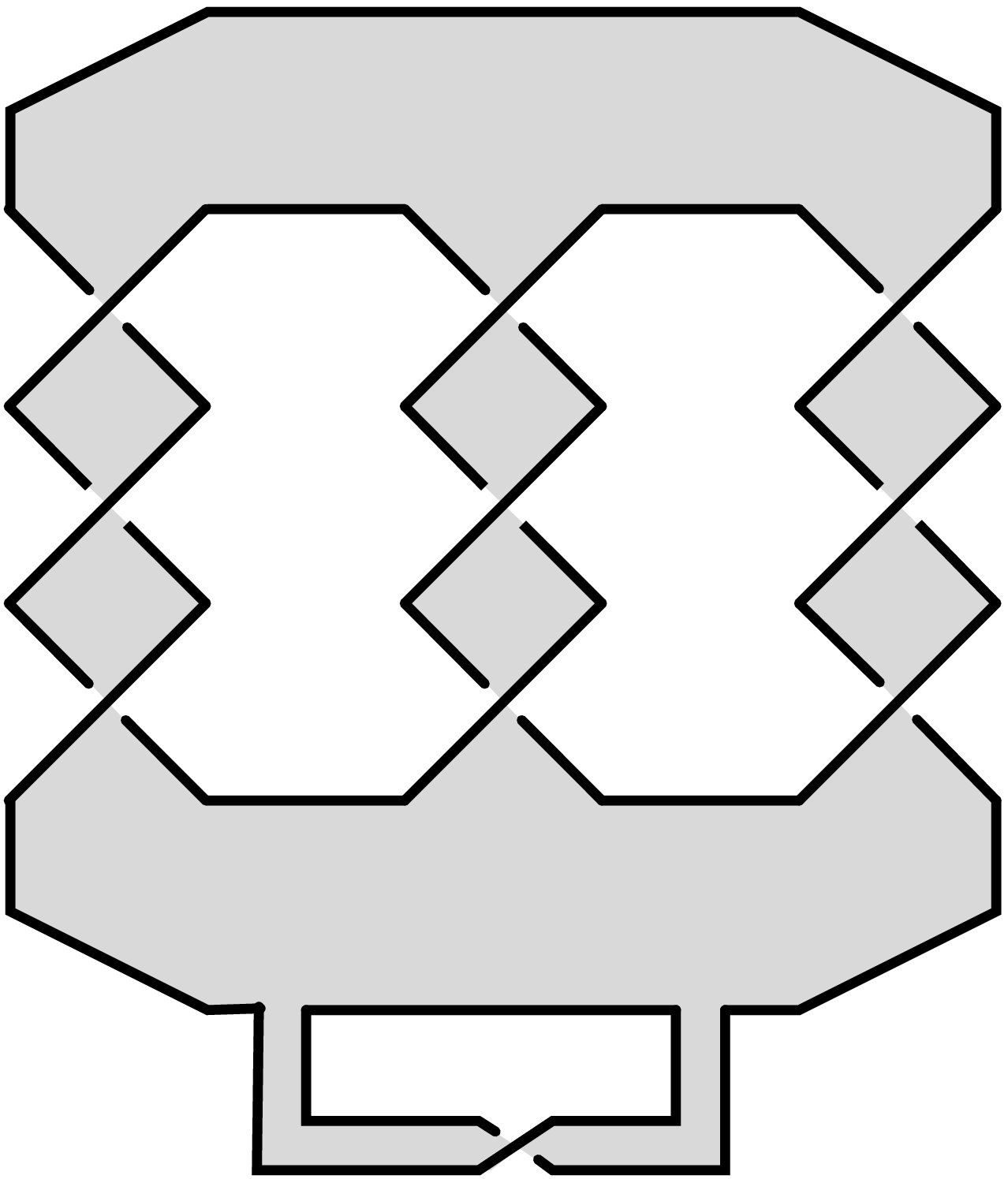}}}
  \end{picture}
 \end{center}
 \caption{$F_A=F_B$ for the $(2,3,3)$-pretzel knot in the class (a)
          and $F_B$ for the $(3,3,3)$-pretzel knot in the class (b)}
\label{Fig:Surfaces}
\end{figure}

%
For a pretzel knot $K=P(p_1,p_2,\ldots,p_\NumTangles)$ satisfying (\ref{Eq:Not01Condition}),
a basic edgepath system $\BasicEdgepathSystem_A=\{\BasicEdgepath_{A,i}\}$
whose edgepaths are given by 
\[
\BasicEdgepath_{A,i}=
\BasicEdgepath_{i,a}=\angleb{0}\textrm{\,--\,}\angleb{1/p_i}
,
\]
is regarded as a type II edgepath system $\EdgepathSystem_A=\{\Edgepath_{A,i}\}$.
Though a candidate edgepath system corresponds to multiple candidate surfaces in general,
by construction, 
all the surfaces obtained from $\EdgepathSystem_A$ are isotopic to each other. 
Let $F_A$ denote the surface so obtained, which is in fact a spanning surface.
$F_A$ is obtained by naturally spanning
a standard diagram of the pretzel knot $K$
as in the left figure of Figure \ref{Fig:Surfaces}.
If the tuple $(p_1,p_2,\ldots,p_\NumTangles)$ corresponding to $K$ 
satisfies the condition (a) in Introduction,
the surface $F_A$ so constructed is non-orientable. 
Meanwhile, if the condition (b) is satisfied, the surface $F_A$ becomes orientable.

We define a non-orientable spanning surface $F_B$ as follows: 
If the condition (a) is satisfied, we set $F_B=F_A$, 
otherwise $F_B$ is obtained by attaching a half-twisted band to $F_A$.
%
%
See Figure \ref{Fig:Surfaces}.
By considering the first Betti number of this $F_B$, we obtain 
\[
\Crosscap(K) \le \Betti ( F_B ) = 
\left\{
\begin{array}{lcl}
\NumTangles-1 & \hspace{1cm} & \mbox{when } (\textrm{a}) \mbox{ is satisfied, }\\ 
\NumTangles & & 				\mbox{when } (\textrm{b}) \mbox{ is satisfied.}
\end{array}
\right.
\]
Hence, Theorem \ref{Thm:Main:Crosscap} claims that the equality in the inequality above holds; 
in other words, $F_B$ attains the crosscap number.

\begin{remark}
If $\NumTangles=2$, then $K$ is $(2,n)$-cabled and $\Crosscap(K)=1$. 
Thus Theorem \ref{Thm:Main:Crosscap} holds in this case. 
If $\NumTangles=1$, then $K$ is the trivial knot, 
and its crosscap number is defined as $\Crosscap(K)=0$.
Thus, the expression of Theorem \ref{Thm:Main:Crosscap} is not true 
when $N=1$ and the condition (b) with $p_1$ odd is satisfied. 
\end{remark}


\section{proof}
\label{Sec:Proof}


%
%
%

As claimed in \cite{IM}, in the algorithm of Hatcher and Oertel, 
the ratio $\frac{-\chi}{\sharp s}  (F)$ of the negative $-\chi (F)$ of $\chi (F)$
and $\sharp s (F)$ appears to be more natural 
for a surface $F$. 
Hence we will perform estimation for $\frac{-\chi}{\sharp s} $, and 
our goal is to give the following.

%
\begin{proposition}
\label{Prop:Bound:Euler}
Assume that
$K$ is a non-trivial pretzel knot $P(p_1,p_2,\ldots,p_\NumTangles)$
with $\NumTangles\ge 3$
satisfying (\ref{Eq:Not01Condition}).
Let $F_B$ be the embedded surface in the exterior of $K$
defined in Section \ref{Sec:Preliminary}.
Then,
for any candidate surface $F_E$ in the exterior of $K$, 
the inequality
\begin{eqnarray}
\label{Eqn:Bound:Euler}
\frac{-\chi}{\sharp s} (F_B)-1 \le \frac{-\chi}{\sharp s} (F_E) 
\end{eqnarray}
holds.
Moreover, if the equality holds,
$F_E$ does not satisfy the condition  ``non-orientable and spanning''.
\end{proposition}

Note that $\frac{-\chi}{\sharp s}(F_B)$ is $\NumTangles-2$ for (a) and $\NumTangles-1$ for (b).
Besides, for a spanning surface $F$,
the first Betti number and the Euler characteristic are related by
$\Betti(F)=-\chi(F)+1=\frac{-\chi}{\sharp s}(F)+1$.
With the above proposition,
Theorem \ref{Thm:Main:Crosscap} can be proved as follows.
%

\begin{proof}[Proof of Theorem \textup{\ref{Thm:Main:Crosscap}}]
Since $F_B$ is a non-orientable spanning surface for K, 
it is clear by definition that $\Crosscap(K)\le \Betti(F_B)$.

Assume for a contradiction that $\Crosscap(K)< \Betti(F_B)$ holds.
Then, for the non-orientable spanning surface $F_C$
which attains the crosscap number,
$\Betti(F_C)=\Crosscap(K)\le \Betti(F_B)-1$ holds.
Remark that we cannot tell whether $F_C$ is essential or not. 

If $F_C$ is essential, then 
we have $\Betti(F_B)-1\le \Betti(F_C)$ 
by Proposition \ref{Prop:Bound:Euler}. 
Together with the assumption above, 
this implies that $\Betti(F_B)-1=\Betti(F_C)$ must hold. 
Then, 
by the condition for the equality in Proposition \ref{Prop:Bound:Euler}, 
$F_C$ does not satisfy ``non-orientable and spanning''. 
Though, this contradicts to the assumption that 
$F_C$ is a non-orientable spanning surface.

Now assume that $F_C$ is inessential. 
As claimed in Lemma 3.1 in \cite{HT},
$F_C$ is in fact incompressible. 
Thus it must be boundary-compressible. 
%
By boundary compressions, 
as argued in the proof of Lemma 3.1 in \cite{HT}, 
an essential surface $F_E$ is obtained. 
Then it satisfies $-\chi (F_E)< -\chi (F_C)$, and 
together with 
$\chi(F_E)\le 0$ obtained in \cite{IM}
and 
$\sharp s (F_C) = 1 \le \sharp s (F_E)$, 
we have $\frac{-\chi}{\sharp s}  (F_E)<\frac{-\chi}{\sharp s}  (F_C)$.
Again, by Proposition \ref{Prop:Bound:Euler},
we have $\frac{-\chi}{\sharp s}  (F_B)-1\le \frac{-\chi}{\sharp s}  (F_E)<\frac{-\chi}{\sharp s}  (F_C)$. 
It then follows that $\Betti (F_B)-1 < \Betti (F_C)$.
Though, it contradicts to $\Betti(F_C)\le \Betti(F_B)-1$.

Therefore we conclude that $\Crosscap(K)=\Betti(F_B)$. 
This completes the proof of the theorem. 
\end{proof}

%
%
%

%
%
%

%

The rest of the paper is devoted to proving Proposition \ref{Prop:Bound:Euler}. 
In the following, we will assume that the readers are rather familiar 
with the work of Hatcher and Oertel. 
Please see \cite{HO} or \cite{IM} for details. 

We here remark that 
we will not check whether a candidate surface 
is actually essential or not. 
That is, we will establish our estimation for all candidate surfaces 
including both essential and inessential ones.

In the following, we divide the argument into two main cases 
concerning types of edgepath systems.

\subsection{Type II and type III edgepath systems}


Let us start with considering type II and type III edgepath systems. 
This case is much easier than the case of type I. 

\begin{lemma}
\label{Lem:Bound:typeII-III}
Assume that 
$K$ is a non-trivial pretzel knot $P(p_1,p_2,\ldots,p_\NumTangles)$ with $\NumTangles\ge 3$ satisfying (\ref{Eq:Not01Condition}).
Let $F_B$ be the embedded surface in the exterior of $K$
defined in Section \ref{Sec:Preliminary}.
Then,
for any embedded surface $F_E$ in the exterior of $K$
corresponding to a type II or III edgepath system,
the inequality
\begin{eqnarray*}
\frac{-\chi}{\sharp s} (F_B)-1 \le \frac{-\chi}{\sharp s} (F_E) 
\end{eqnarray*}
holds.
Moreover, if the equality holds,
$F_E$ does not satisfy the condition  ``non-orientable and spanning''.
\end{lemma}

An edgepath system $\EdgepathSystem$ corresponds to multiple surfaces in general.
Though, the value of $\frac{-\chi}{\sharp s} (F)$ 
is common for any surface $F$ of these surfaces.
See \cite{IM} for example. 
Thus, in the following, 
we denote the common value by $\frac{-\chi}{\sharp s} (\EdgepathSystem)$ for brevity. 


\begin{proof}[Proof of Lemma \textup{\ref{Lem:Bound:typeII-III}}]
Recall that, 
for each edgepath $\Edgepath_i$ in an edgepath system 
$\EdgepathSystem=(\Edgepath_1, \Edgepath_2,\ldots ,\Edgepath_\NumTangles)$, 
the length $|\Edgepath_{i,>0}|$ of 
the part $\Edgepath_{i,>0}$ lying in the region $u>0$ is one or more, 
if $\EdgepathSystem$ is of type II or III. 
On the other hand, 
a formula of calculating $\frac{-\chi}{\sharp s}  (\EdgepathSystem)$ 
for $\EdgepathSystem$ is given as 
\begin{eqnarray*}
\frac{-\chi}{\sharp s} (\EdgepathSystem)
&=&
 \sum_{i=1}^{\NumTangles}|\Edgepath_{i,>0} |
\label{Eq:Formula:EulerCharTypeIII}
\end{eqnarray*}
when $\EdgepathSystem$ is of type III, 
and as 
\begin{eqnarray*}
\frac{-\chi}{\sharp s} (\EdgepathSystem)
&=&
 \left( \sum_{i=1}^{\NumTangles}(|\Edgepath_{i,>0} | \right)
 +V -2
\label{Eq:Formula:EulerCharTypeII}
\end{eqnarray*}
when $\EdgepathSystem$ is of type II. 
Here $V$ denotes the number of vertical edges.

From these, we immediately have 
the inequality (\ref{Eqn:Bound:Euler}) in Proposition \ref{Prop:Bound:Euler}; 
\[
\frac{-\chi}{\sharp s} (\EdgepathSystem) \ge \NumTangles-2
.
\]

The equality can hold only when the condition (b) is satisfied and 
the edgepath system $\EdgepathSystem$ is of type II with no vertical edges 
whose all edgepaths have length 1 in the region $u>0$. 
Note that other than edgepaths of the form 
$\angleb{0}$\,--\,$\angleb{\pm 1/p}$, 
the only edgepath $\Edgepath=\angleb{\pm 1}$\,--\,$\angleb{\pm 1/2}$
satisfies $|\Edgepath|=1$.
However, 
since we here think about the case where the condition (b) is satisfied, 
i.e., all $p_i$'s are odd, 
we do not have to consider $\angleb{\pm 1}$\,--\,$\angleb{\pm 1/2}$ for $\EdgepathSystem$. 
Therefore the equality holds 
only when all edgepaths in $\EdgepathSystem$ 
are of the form $\angleb{0}$\,--\,$\angleb{\pm 1/p}$. 
This $\EdgepathSystem$ corresponds to the surface $F_A$,
which is orientable when the condition (b) is satisfied. 
Hence, the condition for the equality in Proposition \ref{Prop:Bound:Euler} is satisfied.
\end{proof}

%
%
%

\subsection{Type I edgepath systems}


For a type I edgepath system $\EdgepathSystem$,
a formula of calculating $\frac{-\chi}{\sharp s}  (\EdgepathSystem)$ is given as 
\begin{eqnarray}
\frac{-\chi}{\sharp s} (\EdgepathSystem)
&=&
 \sum_{i=1}^{\NumTangles}
  \left(
   \left\{
    \begin{array}{l}
     0 \\
     ~~~~~(\textrm{  if $\Edgepath_i$ is constant }) \\
     |\Edgepath_{i} |\\
     ~~~~~(\textrm{ otherwise }) \\
    \end{array}
   \right.
  \right)
\label{Eq:Formula:EulerCharTypeI}
\\
&& 
+\NumTangles_{\mathrm{const}}-\NumTangles
+\left(
\NumTangles-2-\sum_{\Edgepath_i \in \ConstantEdgepaths}\frac{1}{|p_i|}
\right)
\frac{1}{1-u}
~~.
\nonumber
\end{eqnarray}
Here, $\ConstantEdgepaths$ denotes the constant edgepaths in $\EdgepathSystem$,
$\NumTangles_{\mathrm{const}}$ is the number of the constant edgepaths,
and
$|\Edgepath|$ is the length of an edgepath $\Edgepath$.
The length of an edgepath coincides with the path length in the graph
if all the edges in the edgepath are complete edges.

Originally, this formula should be applied to 
the edgepath systems satisfying the gluing consistency.
Though, we can apply this to other edgepath systems formally.
In this subsection, 
by using the formula, we will show that
$\frac{-\chi}{\sharp s} (\EdgepathSystem)$ is bounded from below for all possibilities
even if we ignore the gluing consistency in most cases.

%
%
%

\subsubsection{Preparations}
We here introduce certain new functions, a variable and a preorder 
to clarify and simply the arguments. 


First let us introduce a function $Y_\BasicEdgepathSystem: (0,1) \longrightarrow \mathbb{R}$ 
for a fixed basic edgepath system $\BasicEdgepathSystem$ in the following way: 
We take
the extended basic edgepath system $\widetilde{\BasicEdgepathSystem}$
for $\BasicEdgepathSystem$.
For $0<u_0<1$,
we cut $\widetilde{\BasicEdgepathSystem}$ at $u=u_0$ 
in the manner described in the previous section,
and let $\EdgepathSystem_{u_0}$ denote the obtained edgepath system,
which is nearly type I but may not satisfy the gluing consistency.
%
Then, we calculate 
the value $\frac{-\chi}{\sharp s}  (\EdgepathSystem_{u_0})$ 
for the edgepath system $\EdgepathSystem_{u_0}$ 
formally by applying the formula (\ref{Eq:Formula:EulerCharTypeI}). 
The function $Y_\BasicEdgepathSystem(u)$ is defined to
take the value $\frac{-\chi}{\sharp s}  (\EdgepathSystem_{u_0})$ for $u=u_0$.

%
%
%


This function $Y_\BasicEdgepathSystem$ becomes much easier to see 
by introducing a new variable $w$ and a new function $x$ as follows: 
We define a new variable $w$ as $w=1/(1-u)$.
By this variable transformation,
we define a new function $X_\BasicEdgepathSystem$ by
$X_\BasicEdgepathSystem(w)=Y_\BasicEdgepathSystem(1-1/w)=Y_\BasicEdgepathSystem(u)$.
The function $x_{\BasicEdgepath}(w) : [1,\infty) \rightarrow \mathbb{R}$ is
defined for a fixed basic edgepath $\BasicEdgepath$ as follows.
We prepare 
the extended basic edgepath $\widetilde{\BasicEdgepath}$ for $\BasicEdgepath$.
For $w>1$,
we obtain an edgepath $\Edgepath_{\BasicEdgepath,w} $ 
by cutting $\widetilde{\BasicEdgepath}$ at $u=1-(1/w)$ in our manner.
Then the function $x_{\BasicEdgepath}$ is defined by; 
\begin{eqnarray*}
x_{\BasicEdgepath}(w)&=&
\left\{
 \begin{array}{ll}
  1-\frac{1}{q} w 
  & (\textrm{  if $\Edgepath_{\BasicEdgepath,w}$ is a constant edgepath} \\
  & \textrm{ \ \ with $v$-coordinate $p/q$}) \\
  \sum_{j=1}^{\textrm{number of edges in $\Edgepath_{\BasicEdgepath,w}$}}
   |\Edge_j| 
  & (\textrm{ if $\Edgepath_{\BasicEdgepath,w}$ is non-constant }) 
 \end{array}
\right.
\end{eqnarray*}
where
\begin{eqnarray*}
|\Edge_j| &=&
   \left\{
    \begin{array}{ll}
     \frac{s_j-w}{s_j-q_j}
       ~~~(= \frac{s_j}{s_j-q_j}
             -\frac{1}{s_j-q_j}w )     
     & \textrm{ ( if $\Edge_j$ is a partial edge of $\angleb{p_j/q_j}$\,--\,$\angleb{r_j/s_j}$} \\
     & \textrm{ \ \ ending at $u=1-1/w$) } \\
     1 
     & \textrm{ ( if $\Edge_j$ is a complete edge) } 
    \end{array}
   \right.
.
\end{eqnarray*}
%
%
Note that 
$x_{\BasicEdgepath}(w)$ is a piecewise linear function 
and is strictly monotonically decreasing.
With these $w$ and $x_{\BasicEdgepath}$, 
the function $X_\BasicEdgepathSystem$ is expressed simply as follows:
\begin{eqnarray*}
X_{\BasicEdgepathSystem}(w)
&=&
 [ (\NumTangles-2)w-\NumTangles ]
+\sum_{i=1}^{\NumTangles} x_{\BasicEdgepath_i} (w),
\end{eqnarray*}
where $\BasicEdgepath_i$'s are the basic edgepaths in 
the basic edgepath system $\BasicEdgepathSystem$. 
%
%
We will always use this expression 
in actual calculations of $X_{\BasicEdgepathSystem}(w)$ in the sequel. 
%
First, we have the following proposition.

%
%
\begin{proposition}
\label{Prop:Order:edgepaths}
For the $(p_1$, $p_2$, $\ldots$, $p_\NumTangles)$-pretzel knot $K$ with $\NumTangles\ge 3$ satisfying (\ref{Eq:Not01Condition}), 
let $\BasicEdgepathSystem_A$ be the basic edgepath system 
$\{\angleb{0}$\,--\,$\angleb{1/p_1}$, 
$\angleb{0}$\,--\,$\angleb{1/p_2}$, $\ldots$, $\angleb{0}$\,--\,$\angleb{1/p_\NumTangles}\}$. 
Then,
for any basic edgepath system $\BasicEdgepathSystem$ for $K$,
the inequality
$
X_{\BasicEdgepathSystem_{A}}(w) \le X_{\BasicEdgepathSystem}(w)$
holds for any $w > 1$. 
\end{proposition}

\begin{proof}
For a rational tangle $1/p_i$, possible basic edgepaths are
$\BasicEdgepath_{i,a}=\angleb{0}$\,--\,$\angleb{1/p_i}$ and
$\BasicEdgepath_{i,b}=\angleb{s_i}$\,--\,$\angleb{s_i/2}$\,--\,$\ldots$\,--\,$\angleb{s_i/(|p_i|-1)}$\,--\,$\angleb{1/p_i)}$,
where $s_i$ denotes the sign of $p_i$.
\begin{figure}[hbt]
 \begin{center}
  \begin{picture}(250,100)
   \put(0,0){\scalebox{1.0}{\includegraphics{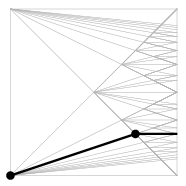}}}
   \put(150,0){\scalebox{1.0}{\includegraphics{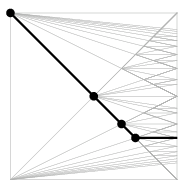}}}
  \end{picture}
 \end{center}
 \caption{$\widetilde{\BasicEdgepath_{i,a}}$ and $\widetilde{\BasicEdgepath_{i,b}}$ for $p_i=4$}
\end{figure}
\\
The functions $x_{\BasicEdgepath}$ corresponding to these extended basic edgepaths are
\begin{eqnarray*}
x_{\BasicEdgepath_{i,a} }(w)
&=&
\left\{
 \begin{array}{lll}
  1 - \frac{1}{|p_i|-1}(w-1) &  & (1 \le w \le |p_i|) \\
  1 - \frac{1}{|p_i|} w & & (|p_i|<w)
 \end{array}
\right.
\end{eqnarray*}
and
\begin{eqnarray*}
x_{\BasicEdgepath_{i,b} }(w)
&=&
\left\{
 \begin{array}{lll}
  |p_i| - w &  & (1 \le w \le |p_i|) \\
  1 - \frac{1}{|p_i|} w & & (|p_i|<w)
 \end{array}
\right.
.
\end{eqnarray*}
\begin{figure}[hbt]
 \begin{center}
  \begin{picture}(328,110)
  \scalebox{0.8}{
   \put(0,0){
    \put(0,10){\scalebox{1.0}{\includegraphics{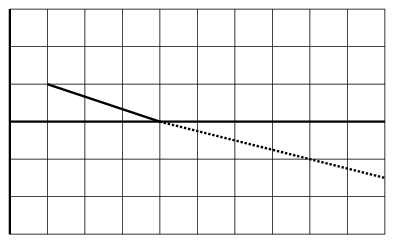}}}
    \put(2,135){$x_{\BasicEdgepath_{i,a}(w)}$}
    \put(-1,125){$3$}
    \put(-1,70){$0$}
    \put(-10,19){$-3$}
    \put(193,77){$w$}
    \put(12,65){$0$}
    \put(191,65){$10$}
   }
  \put(210,0){
   \put(0,10){\scalebox{1.0}{\includegraphics{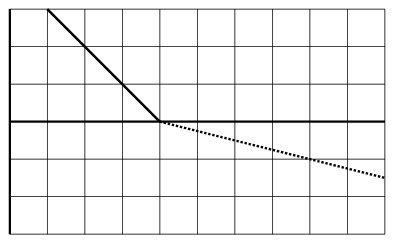}}}
   \put(2,135){$x_{\BasicEdgepath_{i,b}(w)}$}
   \put(-1,125){$3$}
   \put(-1,70){$0$}
   \put(-10,19){$-3$}
   \put(193,77){$w$}
   \put(12,65){$0$}
   \put(191,65){$10$}
  }
  }
  \end{picture}
 \end{center}
 \caption{$x_{\BasicEdgepath_{i,a}}(w)$ and $x_{\BasicEdgepath_{i,b}}(w)$ for $p_i=4$}
\end{figure}
\\
With respect to these two edgepaths
the inequality
$x_{\BasicEdgepath_{i,a}}(w)\le x_{\BasicEdgepath_{i,b}}(w)$ 
clearly holds for any $w>1$. 
Hence, 
for $\BasicEdgepathSystem_A$, 
$X_{\BasicEdgepathSystem_A}(w)$ takes minimum among $X_{\BasicEdgepathSystem}$ of all the basic edgepath systems $\BasicEdgepathSystem$.
\end{proof}

Next we introduce the following preorder for tuples of integers.

%
%
\begin{definition}
Let
$(p_1, p_2, \ldots, p_n)$ and
$(p^\prime_1, p^\prime_2, \ldots, p^\prime_n)$
be $n$-tuples of integers.
Assume that,
for these tuples,
after replacing each entry by its absolute value,
rearrange entries in ascending order,
we obtain $(q_1, q_2, \ldots, q_n)$
and
$(q^\prime_1, q^\prime_2, \ldots, q^\prime_n)$.
Then,
if $q_i\le q^\prime_i$ holds for each $i=1,2,\ldots,n$,
we define $(p_1, p_2, \ldots, p_n) \le (p^\prime_1, p^\prime_2, \ldots, p^\prime_n)$.
Since this relation is reflexive and transitive, this gives a preorder.
Moreover
the binary relation $\sim$ defined as
$x\sim y \Leftrightarrow ( x\le y \textrm{ and } y\le x)$ 
becomes an equivalence relation.
Let $|p_1, p_2, \ldots, p_n|$ denote 
the equivalent class of $(p_1, p_2, \ldots, p_n)$.

For
the $(p_1, p_2, \ldots, p_n)$-pretzel knot $K$
and 
the $(p^\prime_1, p^\prime_2, \ldots, p^\prime_n)$-pretzel knot $K^\prime$,
we define a preorder of the pretzel knots
by $K\le K^\prime \Leftrightarrow (p_1, p_2, \ldots, p_n) \le (p^\prime_1, p^\prime_2, \ldots, p^\prime_n)$.
\end{definition}

With this preorder,
we obtain the following proposition naturally. 

%
%
\begin{proposition}
\label{Prop:Order:knots}
Assume that
 $(p_1$, $p_2$, $\ldots$, $p_\NumTangles)\le (p^\prime_1$, $p^\prime_2$, $\ldots$, $p^\prime_\NumTangles)$.
Let $K$ and $K^\prime$ 
be the $(p_1$, $p_2$, $\ldots$, $p_\NumTangles)$-pretzel knot and 
the $(p^\prime_1$, $p^\prime_2$, $\ldots$, $p^\prime_\NumTangles)$-pretzel knot 
both satisfying (\ref{Eq:Not01Condition})
respectively. 
Let $\BasicEdgepathSystem_A$ denote 
the basic edgepath system 
$\{\angleb{0}$\,--\,$\angleb{1/p_1}$, $\angleb{0}$\,--\,$\angleb{1/p_2}$, $\ldots$ ,$\angleb{0}$\,--\,$\angleb{1/p_\NumTangles}\}$
for the knot $K$,
and
$\BasicEdgepathSystem^\prime$ denote
an arbitrary basic edgepath system for the knot $K^\prime$.

If an inequality
$X_{\BasicEdgepathSystem_{A}}(w) >t$ 
holds for some real number $t$ for any $w > 1$, 
then the inequality
$X_{\BasicEdgepathSystem^\prime}(w) >t$ 
also holds for any $w > 1$. 
The similar fact holds in the case of $X_{\BasicEdgepathSystem_{A}}(w) \ge t$. 


\end{proposition}

\begin{proof}
Let $T$ and $T^\prime$ be the $1/p$ tangle and the $1/p^\prime$ tangle 
satisfying $|p|\le|p^\prime|$. 
For two basic edgepaths 
$\BasicEdgepath_a= \angleb{0}$\,--\,$\angleb{1/p}$ and
$\BasicEdgepath^\prime_a= \angleb{0}$\,--\,$\angleb{1/p^\prime}$
corresponding to $T$ and $T^\prime$ respectively,
we immediately have
$
x_{\BasicEdgepath_a}(w)\le x_{\BasicEdgepath^\prime_a}(w)
$
for any $w>1$. 

This implies that: 
\begin{itemize}
 \item [(1)]
  When $(p_1, p_2, \ldots, p_\NumTangles)\le(p^\prime_1, p^\prime_2, \ldots, p^\prime_\NumTangles)$, we have:
    The inequality 
    $X_{\BasicEdgepathSystem_{A}}(w) \le X_{\BasicEdgepathSystem^\prime_{A}}(w) $ 
    holds for any $w > 1$. 
	Here $\BasicEdgepathSystem^\prime_A$ denotes
	a basic edgepath system $\{\angleb{0}$\,--\,$\angleb{1/p^\prime_1}$, 
	$\angleb{0}$\,--\,$\angleb{1/p^\prime_2}$, $\ldots$ ,
	$\angleb{0}$\,--\,$\angleb{1/p^\prime_\NumTangles}\}$
	for the knot $K^\prime$. 
 \item[(2)]
  When
  $|p_1, p_2, \ldots, p_\NumTangles|
  =|p^\prime_1, p^\prime_2, \ldots, p^\prime_\NumTangles|$, we have:
    The identity
    $X_{\BasicEdgepathSystem_A}(w)=X_{\BasicEdgepathSystem^\prime_A}(w)
    $ holds for any $w > 1$.
\end{itemize}

Together with Proposition \ref{Prop:Order:edgepaths} and 
the definition of our preorder, 
these observations imply the assertion in the proposition. 
\end{proof}

%
%
%


Once we have an appropriate lower bound of the function $X_{\BasicEdgepathSystem_A}$
for the edgepath system $\BasicEdgepathSystem_A$
for a pretzel knot $K$,
by the propositions above,
we also have the same lower bound of $X_{\BasicEdgepathSystem}$
for any edgepath system $\BasicEdgepathSystem$
of any pretzel knot $K^\prime$ equal to or greater than $K$.
For any type I edgepath system 
obtained for these edgepath systems and these pretzel knots,
the value of $\frac{-\chi}{\sharp s}$ is bounded by the same bound naturally.
Thus, we can establish the required lower bound
for all type I edgepath systems of many pretzel knots at the same time,
without solving the equation $\widetilde{\BasicEdgepathSystem}(u)=0$.
We only have to give a case-by-case argument for only some remaining pretzel knots.

\subsubsection{Type I edgepath systems for $\NumTangles \ge 4$}

We first consider the number of tangles $N$ is at least 4. 

%
%

\begin{lemma}
\label{Lem:Bound:typeI:N>3}
Assume that
$K$ is a non-trivial pretzel knot $K=P(p_1,p_2,\ldots,p_N)$ with $N\ge4$
satisfying (\ref{Eq:Not01Condition}).
Let $F_B$ be the embedded surface in the exterior of $K$
defined in Section \ref{Sec:Preliminary}.
Then,
for any embedded surface $F_E$ in the exterior of $K$
corresponding to a type I edgepath system, 
the inequality
\begin{eqnarray*}
\frac{-\chi}{\sharp s} (F_B)-1 \le \frac{-\chi}{\sharp s} (F_E) 
\end{eqnarray*}
holds.
Moreover, if the equality holds,
$F_E$ does not satisfy the condition  ``non-orientable and spanning''.
\end{lemma}

\begin{proof}
Suppose first that $(p_1,p_2,\ldots,p_N)$ satisfies  the condition (a). 
For the pretzel knot $K$ with $(p_1,p_2,\ldots,p_\NumTangles)=(2,3,\ldots,3)$,
and the basic edgepath system $\BasicEdgepathSystem_A = 
\{\angleb{0}$\,--\,$\angleb{1/2}$, $\angleb{0}$\,--\,$\angleb{1/3}$, $\ldots$, $\angleb{0}$\,--\,$\angleb{1/3}\}$,
the function $X_{\BasicEdgepathSystem_A} (w)$ is 
a piecewise linear function described as; 
$\NumTangles-2 (>\NumTangles-3)$ for $w=1$,
$(\NumTangles-5)w/2  + (\NumTangles+1)/2$ for $1\le w\le 2$,
$(3\NumTangles-9)/2 ~~(=\NumTangles-3+(\NumTangles-3)/2>\NumTangles-3)$ for $w=2$,
$(\NumTangles-4)w/2 + (\NumTangles-1)/2$ for $2\le w\le 3$,
$(4\NumTangles-13)/2 ~~(=\NumTangles-3+(2\NumTangles-7)/2>\NumTangles-3)$ for $w=3$,
$(4\NumTangles-13)w/6 $ for $3\le w$,
and
goes to $\infty$ when $w$ goes to $\infty$.
Clearly, $X_{\BasicEdgepathSystem_A}(w)>\NumTangles-3$ holds for any $w>1$.

For any $(p_1,p_2,\ldots,p_\NumTangles)$-pretzel knot in this case,
that is, the conditions (\ref{Eq:Not01Condition}) and (a) are satisfied and $\NumTangles\ge 4$,
$(p_1,p_2,\ldots,p_\NumTangles)$ is equal to or greater than $(2,3,\ldots,3)$. 
Therefore $F_E$ always satisfies
\[
\frac{-\chi}{\sharp s}(F_B)-1=
\NumTangles-3 < \frac{-\chi}{\sharp s}(F_E)
\]
by Proposition \ref{Prop:Order:knots}.

Suppose next that $(p_1,p_2,\ldots,p_N)$ satisfies  the condition (b). 
Since $\NumTangles$ must be odd, we remark that $\NumTangles\ge 5$.
For the pretzel knot $K$ with $(p_1,p_2,\ldots,p_\NumTangles)=(3,3,\ldots,3)$,
and the basic edgepath system $\BasicEdgepathSystem_A = 
\{\angleb{0}$\,--\,$\angleb{1/3}$, $\angleb{0}$\,--\,$\angleb{1/3}$, $\ldots$ ,$\angleb{0}$\,--\,$\angleb{1/3}\}$, 
the function $X_{\BasicEdgepathSystem_A}(w)$ is
a piecewise linear function described as; 
$\NumTangles-2$ for $w=1$,
$(\NumTangles-4)w/2  + \NumTangles/2$ for $1\le w\le 3$,
$2\NumTangles-6 ~~(=\NumTangles-2+(\NumTangles-4)>\NumTangles-2)$ for $w=3$,
$(2\NumTangles-6)w/3 $ for $3\le w$
and
goes to $\infty$ when $w$ goes to $\infty$.
Clearly, $X_{\BasicEdgepathSystem_A}(w)>\NumTangles-2$ holds for any $w>1$. 
%
For any $(p_1,p_2,\ldots,p_\NumTangles)$-pretzel knot in this case,
that is, the conditions (\ref{Eq:Not01Condition}) and (b) are satisfied and $\NumTangles\ge 5$,
$(p_1,p_2,\ldots,p_\NumTangles)$ is equal to or greater than $(3,3,\ldots,3)$. 
Therefore $F_E$ always satisfies
\[
\frac{-\chi}{\sharp s}(F_B)-1=\NumTangles-2 < \frac{-\chi}{\sharp s}(F_E)
\]
by Proposition \ref{Prop:Order:knots}.
\end{proof}



\subsubsection{Type I edgepath systems for $\NumTangles=3$}

The case $N=3$ only remains,
and we treat it in the next lemma. 

%
%
%

%
%

\begin{lemma}
\label{Lem:Bound:typeI:N=3}
Assume that 
$K$ is a non-trivial pretzel knot $K=P(p_1,p_2,p_3)$
satisfying (\ref{Eq:Not01Condition}).
Let $F_B$ be the embedded surface in the exterior of $K$
defined in Section \ref{Sec:Preliminary}.
Then,
for any embedded surface $F_E$ in the exterior of $K$
corresponding to a type I edgepath system, 
the inequality
\begin{eqnarray*}
\frac{-\chi}{\sharp s} (F_B)-1 \le \frac{-\chi}{\sharp s} (F_E) 
\end{eqnarray*}
holds.
Moreover, if the equality holds,
$F_E$ does not satisfy the condition  ``non-orientable and spanning''.
\end{lemma}

\begin{proof}
Suppose first that $(p_1,p_2,p_3)$ satisfies  the condition (a). 
%
If $K$ is $(-2,3,3)$ or $(-2,3,5)$-pretzel knot,
then it is a torus knot.
In this case, as in \cite{IM},
an annulus with $\frac{-\chi}{\sharp s} =0=\NumTangles-3$ is obtained.
However, since the number of sheets $\sharp s=2$ is greater than $1$,
the condition for the equality is satisfied. 
Hence, this case does not matter.
For the other candidate surfaces,
\[
N-3=0<1\le \frac{-\chi}{\sharp s} (F_E)
\]
holds
as in \cite{IM}.
In fact, $\Crosscap(K)=2$ can be obtained also
as the result in (\cite{T04}) about torus knots.
%
If $K$ is one of the other pretzel knots in this case,
by the result given in \cite{IM},
all candidate surfaces satisfy 
\[
\NumTangles-3
=0 <\frac{-\chi}{\sharp s}(F_E)
.
\]

Suppose next that $(p_1,p_2,p_3)$ satisfies  the condition (b). 
%
Similarly in the proof of Lemma \ref{Lem:Bound:typeI:N>3}, 
we prepare $\BasicEdgepathSystem_A$ for 
a $(p_1,p_2,p_3)$-pretzel knot 
satisfying $|p_1,p_2,p_3|=|3,5,7|$ or $|p_1,p_2,p_3|=|5,5,5|$. 
By concrete calculation of the function $X_{\BasicEdgepathSystem_A}$, 
we have the following:
For the pretzel knots with $|3,5,7|$,
the function $X_{\BasicEdgepathSystem_A}(w)$ is 
a piecewise linear function expressed as; 
$1$ for $w=1$,
$w/12 + 11/12$ for $1\le w\le 3$,
$7/6 (>1)$ for $w=3$,
$w/4 + 5/12$ for $3\le w\le 5$,
$5/3 (>1)$ for $w=5$,
$3w/10  + 1/6$ for $5\le w\le 7$,
$34/15 (>1)$ for $w=7$,
$34w/105$ for $7\le w$,
and 
goes to $\infty$ when $w$ goes to $\infty$. 
Similarly, for the pretzel knot with $|5,5,5|$,
the function $X_{\BasicEdgepathSystem_A}(w)$ is expressed as; 
$1$ for $w=1$,
$w/4 + 3/4$ for $1\le w\le 5$,
$2 (>1)$ for $w=5$,
$2w/5$ for $5\le w$,
and goes to $\infty$ when $w$ goes to $\infty$.
The graphs of the functions $X_{\BasicEdgepathSystem_A}(w)$
are illustrated in Figure \ref{Fig:1-3_1-5_1-7}
.
Clearly, $X_{\BasicEdgepathSystem_A}(w)>1$ holds for any $w > 1$. 
Therefore, if $K$ is equal to or greater than either of these two knots with respect to our preorder, 
and if $F_E$ corresponds to a type I edgepath system, 
we have 
\[
\frac{-\chi}{\sharp s}(F_B)-1=\NumTangles-2=1<\frac{-\chi}{\sharp s}(F_E)
\]
by Proposition \ref{Prop:Order:knots}.

\begin{figure}[hbt]

\begin{center}
\begin{picture}(352,136)
\scalebox{0.8}{
\put(0,0){
 \put(0,10){\scalebox{1.0}{\includegraphics{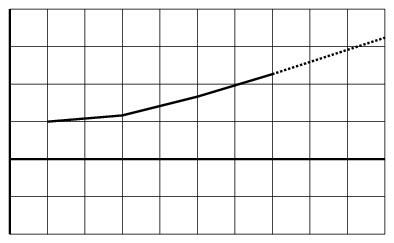}}}
 \put(80,0){1-3,1-5,1-7}
 \put(2,135){$X_{\BasicEdgepathSystem_A}(w)$}
 \put(-1,125){$4$}
 \put(-1,52){$0$}
 \put(-10,19){$-2$}
 \put(193,59){$w$}
 \put(12,47){$0$}
 \put(191,47){$10$}
}
\put(220,0){
 \put(0,10){\scalebox{1.0}{\includegraphics{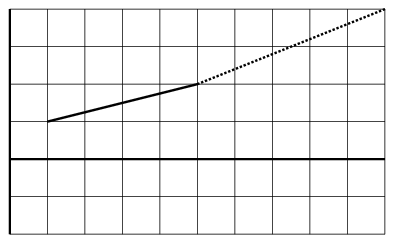}}}
 \put(80,0){1-5,1-5,1-5}
 \put(2,135){$X_{\BasicEdgepathSystem_A}(w)$}
 \put(-1,125){$4$}
 \put(-1,52){$0$}
 \put(-10,19){$-2$}
 \put(193,59){$w$}
 \put(12,47){$0$}
 \put(191,47){$10$}
}
}
\end{picture}
\end{center}

\caption{The graphs of the functions $X_{\BasicEdgepathSystem_A}(w)$ for $\BasicEdgepathSystem_A$ of pretzel knots with $|3,5,7|$ and $|5,5,5|$ }
\label{Fig:1-3_1-3_1-3}
\label{Fig:1-3_1-3_1-11}
\label{Fig:1-3_1-5_1-5}
\label{Fig:1-3_1-5_1-7}
\label{Fig:1-5_1-5_1-5}
\end{figure}

The case of $|p_1,p_2,p_3|=|3,3,n|$ with $n\ge3$ and $|p_1,p_2,p_3|=|3,5,5|$ are only remaining. 
Even for these cases, by \cite{IM},
\[
\frac{-\chi}{\sharp s}(F_B)-1=\NumTangles-2=1\le\frac{-\chi}{\sharp s}(F_E)
\]
holds.
Hence, only the surfaces satisfying the equality in the proposition matter.
We investigate these cases independently as follows. 

Note here that 
for $(p_1, p_2, p_3)$ and 
$(p^\prime_1, p^\prime_2, p^\prime_3)=(-p_1, -p_2, -p_3)$,
there exists a one-to-one correspondence
between candidate surfaces for the former
and candidate surfaces for the latter
preserving the number $\sharp s$ of sheets 
and the Euler characteristic $\chi$.
Similarly, 
if $(p^\prime_1, p^\prime_2, p^\prime_3)$ is obtained
from $(p_1, p_2, p_3)$ by exchanging $p_i$ and $p_j$,
a similar correspondence exists.
Besides, if all $p_i$'s have the common sign,
by the gluing consistency,
there exist no type I edgepath systems.
Eventually,
it suffices to show for the cases of $(p_1, p_2, p_3)=(-3,\pm 3,n)$ and $(-3,\pm5,5)$.
One more important fact here is that if some edgepath in an edgepath system includes a partial edge of a non-horizontal edge, then the number $\sharp s$ of sheets must be greater than $1$,
that is,
the corresponding surface is not spanning.
%


\begin{itemize}

 \item $(p_1, p_2, p_3)=(-3,-3,n)$: \\
  The edgepath system 
  $\EdgepathSystem=\{$
      $((1/2)\angleb{0}+(1/2)\angleb{-1/3})$\,--\,$\angleb{-1/3}$,
      $( (1/2)\angleb{0}+(1/2)\angleb{-1/3})$\,--\,$\angleb{-1/3}$,
      $\angleb{1/2}\textrm{\,--\,}\angleb{1/3}$\,--\,$\ldots$\,--\,$\angleb{1/(n-1)}\textrm{\,--\,}\angleb{1/n}$%
  $\}$
  is the unique type I edgepath system,
  and we have $\frac{-\chi}{\sharp s} (\EdgepathSystem) =n -2 \ge 1$. 
  However, 
  since $\sharp s >1$, the corresponding surface is not a spanning surface.

 \item $(p_1, p_2, p_3)=(-3,3,n)$: \\
  The edgepath system
  $\EdgepathSystem=\{$
      $( (4/(n+1))\angleb{-1/2}+ ((n-3)/(n+1))\angleb{-1/3})$\,--\,$\angleb{-1/3}$,
      $( (2/(n+1))\angleb{0}+ ((n-1)/(n+1))\angleb{1/3})$\,--\,$\angleb{1/3}$,
      $( ((n-1)/(n+1))\angleb{0}+(2/(n+1))\angleb{1/n})$\,--\,$\angleb{1/n}$%
  $\}$
  is the unique type I edgepath system,
  and we have $\frac{-\chi}{\sharp s} (\EdgepathSystem) =1=\NumTangles-2$.
  However, since $\sharp s >1$, 
  the corresponding surface is not a spanning surface.

 \item $(p_1, p_2, p_3)=(-3,-5,5)$: \\
  The edgepath system
  $\EdgepathSystem=\{$
      $( (1/3)\angleb{0}+(2/3)\angleb{-1/3})$\,--\,$\angleb{-1/3}$,
      $( (2/3)\angleb{0}+(1/3)\angleb{-1/5})$\,--\,$\angleb{-1/5}$,
      $( (2/3)\angleb{1/2}+(1/3)\angleb{1/3})$\,--\,$\angleb{1/3}$\,--\,$\angleb{1/4}$\,--\,$\angleb{1/5}$%
  $\}$
  is the unique type I candidate edgepath system,
  and we have $\frac{-\chi}{\sharp s} (\EdgepathSystem)=3>1$. 

 \item $(p_1, p_2, p_3)=(-3,5,5)$: \\
  Only the family of edgepath systems 
  $\EdgepathSystem=\{$
      $( ((2-3u)/(2-2u))\angleb{0}+(u/(2-2u))\angleb{-1/3})$\,--\,$\angleb{-1/3}$,
      $( ((4-5u)/(4-4u))\angleb{0}+(u/(4-4u))\angleb{1/5})$\,--\,$\angleb{1/5}$,
      $( ((4-5u)/(4-4u))\angleb{0}+(u/(4-4u))\angleb{1/5})$\,--\,$\angleb{1/5}$%
  $\}$
  for $0<u<2/3$
  are obtained from this pretzel knot as type I edgepath systems. 
  These edgepath systems correspond to the non-isolated solution of the gluing consistency,
  and we have $\frac{-\chi}{\sharp s} (\EdgepathSystem) =1$
  for any of the edgepath systems.
  $0<1/2\le(4-5u)/(4-4u)<1$ means $\sharp s>1$, and so, 
  any corresponding surface is not a spanning surface. 

\end{itemize}


This completes the proof of the lemma. 
\end{proof}

Consequently, 
Proposition \ref{Prop:Bound:Euler} follows from 
Lemmas \ref{Lem:Bound:typeII-III}, 
\ref{Lem:Bound:typeI:N>3} and 
\ref{Lem:Bound:typeI:N=3}.

%
%

\end{document}